# Le spectre et la torsion analytique des fibrés en droites sur les tores complexes
## Spectrum and analytic torsion of line bundles over complex tori

A. Berthomieu

*Laboratoire Emile Picard. UMR CNRS n° 5580*
*U.F.R. M.I.G. Bat.1R2. Université Paul Sabatier*
*118 route de Narbonne 31062 Toulouse Cédex*
*e-mail:   berthomi@picard.ups-tlse.fr*

*Résumé:*  on calcule le spectre de l'opérateur de Laplace-Dolbeault sur les formes à coefficients dans un fibré en droites $L$ sur un tore complexe plat lorsque $L$ est muni d'une métrique à courbure parallèle. On en déduit alors la torsion analytique de Ray-Singer [RSi] de $L$, généralisant ainsi des résultats de Bost [Bo] pour les fibrés amples et de Ray et Singer [RSi] pour les fibrés plats. On donne de ces derniers une interprétation géométrique en termes de métriques de Quillen.

*Abstract:*  The spectrum of the Laplace-Dolbeault operator for any line bundle with parallel curvature on a flat complex torus is computed. The Ray-Singer analytic torsion [RSi] is then deduced, generalizing thus Bost's result [Bo] for ample line bundles and Ray-Singer's ones [RSi] for flat bundles, of which we a geometric interpretation is given.

**0) Introduction:**

Soit $V$ un espace vectoriel complexe de dimension complexe $n$ muni d'une forme hermitienne définie positive $g$ et $U$ un réseau dans $V$, $H$ une forme hermitienne sur $V$ telle que $E = \operatorname{Im} H$ prenne des valeurs entières sur $U \times U$ et $\alpha : U \to \mathbf{U}(1)$ vérifiant:

$$\alpha(u_1 + u_2) = \alpha(u_1)\,\alpha(u_2)\,\exp i\pi E(u_1, u_2) \tag{1}$$

On considère le tore complexe $T = V/U$ et le fibré $L = L(H, \alpha)$ sur $T$ donné comme quotient de $\mathbf{C} \times V$ par l'action de $U$:

$$\Phi_u(\lambda, z) = \left(\lambda\,\alpha(u) \exp\left[\pi H(z, u) + \frac{\pi}{2} H(u, u)\right], z + u\right) \tag{2}$$

Le théorème d'Appell-Humbert (cf. [Mm page 20]) garantit que tous les fibrés en droites sur $T$ sont obtenus ainsi. On munit $T$ de la métrique parallèle $g$ et le fibré $L$ de la métrique hermitienne naturelle:

$$(\lambda_1, \lambda_2) \to h(\lambda_1, \lambda_2) = \lambda_1 \overline{\lambda}_2 \exp[-\pi H(z, z)] \tag{3}$$

Dans les trois premiers paragraphes, on calcule successivement le spectre du Laplacien de Dolbeault sur les formes de type $(0, k)$ à valeurs dans $L$ et la torsion analytique de Ray-Singer $T_0(T, L)$ dans le cas où $H$ a un noyau nul; on retrouve ainsi le résultat de Bost [Bo proposition 4.2] pour les fibrés amples. Ensuite dans les paragraphes 4 et 5, le cas où $H$ a un noyau non nul est traité; on trouve pour la torsion analytique des résultats analogues à [RSi théorèmes 4.1 et 5.2] qui sont ensuite interprétés en termes de métriques de Quillen au paragraphe 6.

**1) Calcul du Laplacien:**

La connexion de Chern, la courbure et le $c_1$ associés à la métrique (3) sont donnés par:
$$\nabla^L = \mathrm{d} - \pi H(\bullet, z) \qquad \left(\nabla^L\right)^2 = 2\pi i E \qquad c_1(L) = -E$$
On pose:
$$H(z_1, z_2) = g([g^{-1}H]z_1, z_2) = g(z_1, [g^{-1}H]z_2)$$
et on diagonalise $g^{-1}H$ dans une base $(e_i)_{1 \leq i \leq n}$ orthonormée pour $g$ en supposant que les valeurs propres (réelles) $(\mu_i)_{1 \leq i \leq n}$ de $g^{-1}H$ sont rangées dans l'ordre croissant. On appelle enfin $e_i^{(1,0)}$ et $e_i^{(0,1)}$ les projections sur $T^{(1,0)}T$ et $T^{(0,1)}T$ des vecteurs $e_i \in T_{\mathbf{R}}T \cong V$. Dans la trivialisation (2), $\overline{\partial}_L$ reste l'opérateur $\overline{\partial}$ habituel et:
$$\overline{\partial}_L^* = -2 \sum_i \nabla^L_{e_i^{(1,0)}} \iota_{e_i^{(0,1)}} = \overline{\partial}_{\text{habituel}}^* + 2\pi \iota_{g^{-1}Hz^{(0,1)}}$$

où $\iota_v$ désigne le produit intérieur par $v \in T^{(0,1)}T$, $z$ désigne la section tautologique de $\mathcal{C}^\infty(V, V)$ et $z^{(0,1)} \in \mathcal{C}^\infty(V, T^{(0,1)}T)$ sa projection sur $T^{(0,1)}T$. Soit $\Delta^{\text{Riem.}}_{\text{habituel}}$ le Laplacien Riemannien habituel sur $V$ muni de la métrique $g$, on déduit du calcul précédent le Laplacien sur les formes de type $(0, k)$:
$$\Box = \overline{\partial}_L \overline{\partial}_L^* + \overline{\partial}_L^* \overline{\partial}_L = \frac{1}{2}\Delta^{\text{Riem.}}_{\text{habituel}} + 2\pi \nabla''_{g^{-1}Hz^{(0,1)}} + 2\pi \iota_{g^{-1}H} = \overline{\Delta} + 2\pi \iota_{g^{-1}H} \quad (4)$$

où $\nabla''$ coïncide avec l'opérateur $\overline{\partial}$ habituel sur $V$ et:
$$\iota_{g^{-1}H} = \sum_{i=1}^n \mu_i\, e^i \wedge \iota_{e_i^{(0,1)}}$$

où $(e^i)_{1 \leq i \leq n}$ désigne la base de $T^{*(0,1)}T$ duale de la base $(e_i^{(0,1)})_{1 \leq i \leq n}$ de $T^{(0,1)}T$. (cf. [Mr formule (2.8)] pour un calcul analogue).

**2) Calcul du spectre:**

De (4), on déduit que les valeurs propres de $\Box$ sont sommes des valeurs propres de l'opérateur scalaire $\overline{\Delta}$ sur $\mathcal{C}^\infty(V, \mathbf{C})/\Phi$ et de celles de $2\pi \iota_{g^{-1}H}$ sur $\wedge^\bullet T^{*(0,1)}T$.

On suppose que les $\mu_i$ sont linéairement indépendantes sur $\mathbf{Q}$ (en particulier aucune ne s'annule). La diagonalisation de $\iota_{g^{-1}H}$ est donnée par:
$$\iota_{g^{-1}H}(e^{i_1} \wedge e^{i_2} \wedge \cdots \wedge e^{i_k}) = \left(\sum_{j=1}^k \mu_{i_j}\right)(e^{i_1} \wedge e^{i_2} \wedge \cdots \wedge e^{i_k}) \quad (5)$$

pour tous $\{i_1,\cdots,i_k\} \subset \{1,2,\cdots,n\}$. On suppose que $\mu_p$ est la plus grande valeur propre strictement négative de $g^{-1}H$. (La signature de $H$ est donc $n-p$).

**Théorème 1:**
$$\mathrm{Spec}\overline{\Delta} = \left\{2\pi \sum_{i=1}^n n_i|\mu_i| \bigg/ n_i \in \mathbf{N} \text{ et } n \neq 0 \text{ si } i \leq p\right\}$$

*les valeurs propres ayant toutes la multiplicité commune* $(-1)^p\chi(L)$ *(où $\chi(L)$ désigne la caractéristique d'Euler de $L$).*

Ceci généralise [H théorème 1]

**Remarque:** On reconnaît (à la multiplicité près) le spectre de l'oscillateur harmonique sur $\mathbf{R}^n$:
$$\left(-\frac{1}{2}\Delta^{\text{Euclidien}} + 2\pi^2 \sum_{i=1}^n \mu_i^2 x_i^2\right) - \pi \mathrm{Tr}(g^{-1}H)$$

Une telle analogie est explicitée dans certains cas particuliers (cf. [H] et [Be]).

**Corollaire:** *on en déduit le spectre de $\square$:*
$$\mathrm{Spec}\,\square = \left\{2\pi \sum_{i=1}^n n_i|\mu_i| \bigg/ n_i \in \mathbf{N} \text{ pour tout i}\right\}$$

*et la dimension de l'espace propre $E_\lambda^k$ associé à $\lambda = 2\pi \sum_{i=1}^n n_i|\mu_i|$ sur les formes de degré $(0,k)$ est donnée par:*

$$\dim E_\lambda^k = |\chi(L)| C_{\#\{i/n_i \neq 0\}}^{k-\#\{i \leq p/n_i=0\}} \tag{6}$$

*(où le coefficient binomial $C_n^p$ est pris nul si $p<0$ ou $p>n$)*

**Preuve du corollaire:** D'après la décomposition spectrale (5):

$$E_\lambda^k = \bigoplus_{i_1 \leq \cdots \leq i_k} \left[\left(\bigwedge_{j=1}^k e^{i_j}\right) \otimes E_{\lambda-2\pi\sum_{j=1}^k \mu_{i_j}}^0\right] \tag{7}$$

Pour tous les $\lambda' = 2\pi \sum_{i=1}^n m_i|\mu_i|$ où $m_i \geq 1$ pour tout $i \leq p$, les $E_{\lambda'}^0$ sont tous de même dimension $|\chi(L)|$. Cette dernière condition entraîne que si $i \leq p$ et $n_i = 0$, $e^i$ doit obligatoirement intervenir dans (7), alors que si $i \geq p+1$ et $n_i = 0$, $e^i$ ne peut pas intervenir dans (7). (6) en découle immédiatement.

**Preuve du Théorème 1:** La caractéristique d'Euler:

$$\chi(L) = \sum_{k=0}^n (-1)^k \dim E_0^k = \int_T \frac{(-E)^n}{n!} = \left(\prod_{i=1}^n \mu_i\right) \mathrm{Vol}_g(T)$$

où $\text{Vol}_g(T)$ désigne le volume de $T$ pour la métrique $g$, ne s'annule pas d'après l'hypothèse faite sur les $\mu_i$. La formule (7) nous indique alors que toutes les formes harmoniques sont de degré $p$ et que la plus petite valeur propre correspondante de l'opérateur $\overline{\Delta}$ est $2\pi \sum_{i\leq p} |\mu_i|$.

On utilise ensuite le fait que pour tout $\lambda > 0$ le complexe:

$$0 \to E_\lambda^0 \xrightarrow{\overline{\partial}} E_\lambda^1 \xrightarrow{\overline{\partial}} \cdots \xrightarrow{\overline{\partial}} E_\lambda^n \to 0$$

est acyclique, en particulier:

$$\sum_{k=0}^n (-1)^k \dim E_\lambda^k = 0 \tag{8}$$

Soit $\lambda = 2\pi \sum_{i=1}^n m_i |\mu_i|$ (avec $m_i \geq 1$ pour $i \leq p$). On suppose qu'on a établi le théorème sur l'intervalle $[0, \lambda]$. Soit $\lambda'$ le plus petit réel strictement supérieur à $\lambda$ s'écrivant $2\pi \sum_{i=1}^n n_i |\mu_i|$ (avec $n_i \geq 1$ pour $i \leq p$). Aucun réel $\lambda'' \in ]\lambda, \lambda'[$ ne peut être valeur propre de $\overline{\Delta}$ sinon lui correspondrait en vertu de (7) un espace propre:

$$E^p_{\lambda'' + 2\pi \sum_{i\leq p} \mu_i} = \left( \bigwedge_{i\leq p} e^i \right) \otimes E^0_{\lambda''}$$

avec $\lambda'' + 2\pi \sum_{i\leq p} \mu_i > 0$. Aucune valeur propre de $\overline{\Delta}$ strictement supérieure à $\lambda''$ ne peut produire cette même valeur propre en aucun degré pour $\square$. Il en est de même pour les valeurs propres inférieures ou égales à $\lambda$ d'après l'hypothèse de récurrence et l'hypothèse d'indépendance sur $\mathbf{Q}$ des $\mu_i$. Cette valeur propre strictement positive aurait donc un unique espace propre en degré $p$ ce qui contredit (8).

Enfin soit $\mu = \lambda' + \sum_{i\leq p} \mu_i > 0$. L'hypothèse de récurrence et (7) permettent alors de calculer $\dim E_\mu^k$ (6) pour tout $k \neq p$. Ensuite du fait que pour tout $q \geq 1$:

$$\sum_{k=0}^q (-1)^k C_q^k = 0$$

on déduit de (8) la dimension de $E_\mu^p$. (7) et l'hypothèse de récurrence imposent alors la condition:

$$\dim E_{\lambda'}^0 = |\chi(L)|$$

**Remarque:** Le résultat est identique si on relache l'hypothèse d'indépendance sur $\mathbf{Q}$ des $\mu_i$, en demandant tout de même qu'aucune d'entre elles ne s'annule. Il peut y avoir certaines combinaisons des $\mu_i$ qui coïncident auquel cas les espaces propres correspondants se superposent simplement. La compensation nécessaire à la récurrence intervenant toujours en degré $p$, ces compensations ne peuvent se neutraliser.

**Remarque 2:** Le spectre est indépendant du choix de $\alpha$ vérifiant (1), c'est à dire qu'il ne change pas si on tensorise $L$ par un fibré plat: en effet cette opération revient simplement à translater le fibré $L(H, \alpha)$ si $H$ est sans noyau sur $V$.

**Remarque 3:** Le spectre de $\Box$ sur les formes de type (p,q) est le même que sur les formes de type (0,q) si ce n'est que la multiplicité de toutes les valeurs propres est simplement multipliée par $\dim(\wedge^p T^{*(1,0)}T) = C_n^p$.

### 3) Calcul de la torsion analytique:

On suppose toujours $\mathrm{Ker}H = \{0\}$. Soit $N$ l'opérateur qui multiplie par $k$ les formes de degré $k$ et $P^\perp$ le projecteur orthogonal sur le supplémentaire orthogonal de $\mathrm{Ker}\Box$. On rappelle que la torsion analytique de Ray et Singer [RSi] est donnée par:

$$T_0(T, L) = \exp \frac{1}{2}\zeta'_\Box(0)$$

après prolongement méromorphe convenable à $\mathbf{C}$ tout entier de la fonction

$$\zeta_\Box(s) = \mathrm{Tr}\big[(-1)^N N \, \Box^{-s} P^\perp\big] = \sum_{\lambda \in \mathrm{Spec}\Box \setminus \{0\}} \lambda^{-s} \left(\sum_{k=0}^n (-1)^k k \dim E_\lambda^k\right)$$

Or on sait que si $q \geq 2$

$$\sum_{k=0}^q (-1)^k k \, C_q^k = 0 \tag{9}$$

et donc en vertu de (6) et (7), seules interviennent dans le calcul de la torsion analytique les valeurs propres de $\Box$ de type $2\pi n |\mu_i|$ avec $n \geq 1$:

$$\zeta_\Box(s) = (-1)^p |\chi(L)| \, (2\pi)^{-s} \zeta(s) \left(\sum_{i \leq p} |\mu_i|^{-s} - \sum_{i \geq p+1} |\mu_i|^{-s}\right)$$

où $\zeta$ désigne la fonction $\zeta$ de Riemann. On en déduit:

$$\begin{aligned}
\zeta'_\Box(0) &= \chi(L) \left[(2p-n)\zeta'(0) - \zeta(0) \mathrm{Log}\left(\frac{\prod_{i \leq p} 2\pi|\mu_i|}{\prod_{i \geq p+1} 2\pi|\mu_i|}\right)\right] \\
&= -\frac{1}{2} \chi(L) \mathrm{Log}\left(\prod_{i=1}^n |\mu_i|^{\mathrm{sgn}\mu_i}\right)
\end{aligned}$$

et le

**Théorème:**

$$T_0(T, L) = \left(\prod_{i=1}^n |\mu_i|^{\mathrm{sgn}\mu_i}\right)^{-\frac{1}{4}\chi(L)}$$

En particulier si $L$ est ample (i.e. $p = 0$), on retrouve la formule de Bost [Bo prop.4.2]:

$$T_{\text{Bost}}(T, L) = -2 \operatorname{Log} T_0(T, L) = \frac{1}{2} \chi(L) \operatorname{Log}\left(\frac{\chi(L)}{\operatorname{Vol}_g(T)}\right)$$

La première égalité est due à des conventions différentes, $\chi(L)$ ici est égal à $\rho(L)$ dans [Bo], et le facteur $(2\pi)^g$ dans [Bo] correspond à la différence de convention pour la forme de Kähler: $\omega_{\text{Bost}} = (2\pi)^{-1}\omega_g$.

**Remarque:** D'après la remarque 3 on obtient pour tout $p$:

$$T_p(T, L) = T_0(T, L)^{C_n^p} = \left(\prod_{i=1}^{n} |\mu_i|^{\operatorname{sgn}\mu_i}\right)^{-\frac{1}{4}\chi(L)C_n^p}$$

### 4) Le cas où $H$ a un noyau non nul:

On pose $V' = \operatorname{Ker} H \subset V$ et $V'' = V/V'$. $H$ induit naturellement une forme hermitienne sur $V''$ qu'on notera toujours $H$.

**Lemme:** $U' = U \cap V'$ est un réseau sur $V'$. (cf.[W §VI lemmes 2 et 3])

**Preuve:** Soit $U''$ l'image par la projection naturelle $\pi'' : V \to V''$ du réseau $U$. Pour tous $u_1, u_2 \in U$:

$$E(\pi''(u_1), \pi''(u_2)) = E(u_1, u_2)$$

donc $E$ prend des valeurs entières sur $U''$. Comme $E$ est une 2-forme antisymétrique sans noyau sur $V''$, $U''$ est discret et c'est donc un réseau dans $V''$. On en déduit facilement le lemme.

Ce lemme nous enseigne que $T$ est une fibration sur le tore $T'' = V''/U''$ de fibre modelée sur $T' = V'/U'$. On déduit de (2) que la restriction de $L = L(H, \alpha)$ à toutes les fibres de $T$ donne le même fibré plat $P$ sur $T'$ associé au cocycle $\alpha|_{U'}$. Soit $\widehat{V}'$ l'antidual de $V'$ et $\widehat{U}'$ le réseau dans $\widehat{V}'$ dual de $U'$:

$$\widehat{V}' = \operatorname{Hom}_{\mathbf{C}-\text{antilinéaires}}(V', \mathbf{C})$$
$$\widehat{U}' = \{\ell \in \widehat{V}' \,/\, \operatorname{Im}\ell(u) \in \mathbf{Z} \text{ pour tout } u \in U'\}$$

tous les $\ell \in \widehat{V}'$ tels que pour tout $u \in U'$ on ait:

$$\exp(2\pi i \operatorname{Im}\ell(u)) = \alpha(u) \qquad (10)$$

se projettent sur le même point du tore $\widehat{T}' = \widehat{V}'/\widehat{U}'$ dual de $T'$ ($\widehat{T}'$ est la Jacobienne de $T'$, ce point de $\widehat{T}'$ est la classe d'isomorphisme du fibré plat $P$ sur $T'$). On fixe un tel $\ell$ qu'on appelle $\ell_\alpha$ et on prolonge $\alpha$ à $V'$ par:

$$\alpha(z) = \exp(2\pi i \operatorname{Im}\ell_\alpha(z))$$

**Théorème:** *(spectre des fibrés plats)* La décomposition spectrale du Laplacien sur $P$ est alors donnée par la décomposition en série de Fourier:

$$\begin{aligned}
\operatorname{Spec} \Box_P &= \{2\pi^2 \|\ell + \ell_\alpha\|^2 \, / \, \ell \in \widehat{U}'\} \\
E^\bullet_{2\pi^2\|\ell+\ell_\alpha\|^2} &= (\wedge^\bullet T^*T') \otimes \mathbf{C}\gamma_\ell \alpha \\
\text{où} \qquad \gamma_\ell(z) &= \exp(2\pi i \operatorname{Im}\ell(z))
\end{aligned} \qquad (11)$$

En effet, dans une trivialisation de type (2) de $P$, la métrique sur $P$ est triviale et $\Box_P$ est simplement le Laplacien de Dolbeault habituel.

Soit $L'' = L(H, \beta)$ un fibré en droites sur $T''$ associé à $H$ et à un $\beta$ quelconque vérifiant une condition de type (1), et soit $\widetilde{\Box}$ le Laplacien de Dolbeault sur le produit $T' \times T''$ pour le fibré $\pi_1^* P \otimes \pi_2^* L''$ muni de sa métrique naturelle (3), $\pi_1$ et $\pi_2$ étant les projections évidentes, $V'$ et $V''$ étant munis des métriques $g|_{V'}$ et $g|_{V'^\perp}$,

**Théorème 4:** Les spectres de $\Box$ et $\widetilde{\Box}$ coïncident.

**Corollaire:** La torsion analytique est donnée par [RSi théorème 3.3]:

$$T_0(T, L) = T_0(T', P)^{\chi(L'')}$$

et donc elle vaut 1 si $\dim T \geq 2$ [RSi théorème 5.2] (à cause d'une symétrie évidente cf. (9) et (11)), et sinon son logarithme est au coefficient multiplicatif $\chi(L'')$ près donné par une fonction $\theta$ calculée par Ray et Singer [RSi théorème 4.1] (cf. §6).

**5) Preuve du théorème 4:**

Une fois $V''$ muni de $g|_{V^\perp}$, la propriété à démontrer est évidente pour $\iota_{g^{-1}H}$. Il suffit donc de la vérifier en degré 0, donc pour les sections de $L$. Les métriques de type (3) sur les fibrés en droites sont compatibles aux opérations d'image inverse et de produit tensoriel; l'opérateur $\overline{\Delta}$ respecte les directions orthogonales, la difficulté provient de ce que le cocycle (2) ne les respecte pas. On utilise alors une technique analogue à celle de [St §3].

Soit $\ell \in \widehat{U}$ et $\pi'$ la projection orthogonale de $V$ sur $V'$. On pose pour $z \in V$:

$$\delta_\ell(z) = \gamma_\ell(\pi'z)\alpha(\pi'z)$$

Cette fonction est constante dans les directions de $V'^\perp$. Elle permet de définir un fibré plat $P_\ell$ sur $T$ dont la restriction à toutes les fibres de $T$ coïncide avec $P$ via l'action de $U$ sur $\mathbf{C} \times V$:

$$\varphi_u(\lambda, z) = (\lambda \delta_\ell(u), z + u)$$

On considère la fonction sur $U''$:

$$\beta_\ell([u]) = \alpha(u).\delta_\ell(u)^{-1}$$

qui vérifie une condition de type (1). $\big(\beta_\ell([u])$ ne dépend pas du représentant $u \in U$ choisi de $[u] \in U''\big)$. Le fibré $L''_\ell = L(H, \beta_\ell)$ sur $T''$ vérifie alors:

$$\pi''^* L''_\ell \otimes P_\ell \cong L$$

A toute section propre $s$ de $L''_\ell$ associée à une valeur propre $\lambda$ du Laplacien sur $T''$ on associe la section $\pi''^* s \otimes \delta_\ell$ de $L$ sur $V'$, qui est alors une section propre de $\square$ associée à $\lambda + 2\pi^2 \|\ell + \ell_\alpha\|^2$.

Réciproquement, toute section de $L$ se décompose en somme $L^2$ de tels objets: en effet, soit $\sigma$ une section de $L$, alors dans la trivialisation (2), $\sigma \delta_0^{-1}$ est une fonction périodique par rapport à $U'$. On la décompose en série de Fourier dans les directions de $V'$:

$$\sigma = \sum_{\ell \in \widehat{U}} s_\ell \delta_\ell$$

$s_\ell$ est une fonction sur $V'^\perp$ dont la projection sur $V''$ définit une section de $L''_\ell$, que l'on peut alors décomposer selon les sections propres de $L''_\ell$. Enfin, comme tous les $L''_\ell$ sont isospectraux (cf. remarque 2), le théorème est démontré.

**Remarque:** On est ici dans une situation où la formule du produit [RSi théorème 3.3] s'applique. Un cas très général de comparaison de torsions analytiques pour les fibrations est effectué dans [BeBi].

### 6) Une interprétation de la torsion analytique sur les fibrés plats:

On suppose ici que $H = 0$ et donc que $V = V'$ et $U = U'$. Sur le produit des tores mutuellement duaux $T \times \widehat{T}$, le fibré en droites de Poincaré $P_T$ est le quotient de $\mathbf{C} \times V \oplus \widehat{V}$ par l'action de $U \times \widehat{U}$:

$$\phi_{u_1, \hat{u}_2}(\lambda, z_1, \hat{z}_2) = \Big(\lambda \exp \pi \big(\hat{z}_2(u_1) + \overline{\hat{u}_2(z_1)} + \overline{\hat{u}_2(u_1)}\big), z_1 + u_1, \hat{z}_2 + \hat{u}_2\Big)$$

(cf. [Mm page 86]). Il est uniquement caractérisé à isomorphisme près par:

(i) $P_T|_{T \times \{\ell\}}$ est plat et appartient à la classe d'isomorphisme définie par $\ell$.

(ii) $P_T|_{\{0\} \times \widehat{T}}$ est trivial sur $\widehat{T}$.

On considère la projection triviale $\hat{\pi} : T \times \widehat{T} \to \widehat{T}$. Les torsions analytiques de tous les fibrés plats $P_T|_{T \times \{\ell\}}$ sur $T$ interviennent dans l'expression de la métrique de Quillen sur le fibré $\lambda = \det^{-1}(R^\bullet \hat{\pi}_* P_T)$ sur $\widehat{T}$: $\lambda_\ell \cong \otimes_i (\det H^i(T, P_T|_{T \times \{\ell\}}))^{(-1)^{i+1}}$ est muni en tout point $\ell \in \widehat{T}$ d'une norme dite "norme $L^2$" provenant du produit scalaire $L^2$ sur $C^\infty(T, \wedge^\bullet T^{*(0,1)}T \otimes P_T|_{T \times \{\ell\}})$ restreint aux formes harmoniques. La norme de Quillen sur $\lambda_\ell$ est alors donnée par:

$$\| \ \|_{\text{Quillen}} = T^0(T, P_T|_{T \times \{\ell\}}).| \ |_{L^2}$$

Pour tout $\ell \neq 0 \in \widehat{T}$, le fibré plat $P_T|_{T\times\{\ell\}}$ est non trivial, donc sa cohomologie est nulle en tous degrés. Il s'ensuit que le faisceau $R^\bullet \hat{\pi}_* P_T$ a son support réduit à $\{0\} \subset \widehat{T}$. En particulier, $\lambda$ est trivialisé sur $\widehat{T}\backslash\{0\}$, il y admet une section holomorphe $\sigma$ partout non nulle de norme $L^2$ constante égale à 1.

Si dim$T \geq 2$, on déduit du théorème de Hartogs que $\lambda$ est le fibré trivial sur $\widehat{T}$ et $\sigma$ une section triviale. Une application simple du théorème de courbure de Bismut, Gillet et Soulé [BiGSo théorème 0.1] montre que la courbure de la métrique de Quillen sur $\lambda$ est nulle. Cette métrique est donc triviale. La norme $L^2$ l'est aussi, donc la torsion analytique est constante (indépendante de $\ell \in \widehat{T}$). [RSi théorème 5.2] nous enseigne que cette constante est 1.

Si dim$T = 1$, la même application du théorème de courbure [BiGSo] nous montre que le $c_1$ de la métrique de Quillen sur $\lambda$ est égal à la forme volume parallèle sur $\widehat{T}$ de volume total 1. $\lambda$ est donc le fibré en droites sur $\widehat{T}$ associé au diviseur effectif $[0]$, et $\sigma$ prolongée par $\sigma(0) = 0$ est une section globale de $\lambda$: si $T = \mathbf{C}/\Gamma$ où $\Gamma$ est le réseau engendré par 1 et $\tau$ avec Im$\tau > 0$, alors $\widehat{T} = \mathbf{C}/\widehat{\Gamma}$ où $\widehat{\Gamma} = \frac{1}{\mathrm{Im}\tau}\Gamma$; la trivialisation (2) de $\lambda$ et la section $\sigma$ sont données par:

$$H(\hat{z}_1, \hat{z}_2) = \hat{z}_1\overline{\hat{z}_2}\mathrm{Im}\tau \qquad \alpha\Big(\frac{m+n\tau}{\mathrm{Im}\tau}\Big) = (-1)^{mn+m+n} \qquad \sigma(\hat{z}) = \gamma\theta_\tau(\hat{z})\, e^{\frac{\pi}{2}(\mathrm{Im}\tau)\hat{z}^2}$$

où $\gamma$ est une constante et:

$$\hat{z} = \frac{-1}{\mathrm{Im}\tau}(u - v\tau) \qquad\qquad \eta(\tau) = e^{\pi i \frac{\tau}{12}}\prod_{k=1}^{\infty}\Big(1 - e^{2\pi i k\tau}\Big)$$

$$\theta_\tau(\hat{z}) = -\eta(\tau)e^{\pi i(\frac{\tau}{6} - \hat{z}\mathrm{Im}\tau)}\prod_{k\in\mathbf{Z}}\Big(1 - e^{2\pi i(|k|\tau + \hat{z}(\mathrm{Im}\tau)\mathrm{sgn}(k+\frac{1}{2}))}\Big)$$

Le $c_1$ étant parallèle, la métrique de Quillen sur $\lambda$ est multiple de celle donnée par la formule (3). Comme $|\hat{z}|_{L^2}$ vaut 1, on a pour tout $\ell \neq 0$:

$$\|\sigma(\hat{z})\|_{\mathrm{Quillen}} = T_0(T, P_T|_{T\times\{\hat{z}\}}) = \Big|\gamma\,\theta_\tau(\hat{z})\,e^{-\frac{\pi}{2}|\hat{z}|^2\mathrm{Im}\tau}\Big|$$

on retrouve à une constante près la formule [RSi théorème 4.1]:

$$T_0(T, P_T|_{T\times\{\hat{z}\}}), = \Big|\frac{\theta_\tau(\hat{z})}{\eta(\tau)}e^{\pi i(\mathrm{Im}\hat{z})^2\tau}\Big|$$

Celle-ci nous indique alors que:

$$|\gamma| = \frac{1}{|\eta(\tau)|}$$

**Références:**